\ifx\shlhetal\undefinedcontrolsequenc\let\shlhetal\relax\fi

\documentstyle[11pt]{article}

\input amssym.def
\input amssym.tex

\newcount\skewfactor
\def\mathunderaccent#1#2 {\let\theaccent#1\skewfactor#2
\mathpalette\putaccentunder}
\def\putaccentunder#1#2{\oalign{$#1#2$\crcr\hidewidth
\vbox to.2ex{\hbox{$#1\skew\skewfactor\theaccent{}$}\vss}\hidewidth}}
\def\name{\mathunderaccent\tilde-3 }

\newcommand{\V}{{\bf V}}

\newcommand{\rest}{{\restriction}}

\newcommand{\can}{{{}^{\textstyle \omega}2}}
\newcommand{\kom}{{}^{\textstyle \omega}\kappa}
\newcommand{\lom}{{}^{\textstyle \omega}\lambda}
\newcommand{\lh}{\ell{\rm g}}

\newcommand{\con}{2^{\aleph_0}}

\newcommand{\rng}{{\rm rang}}
\newcommand{\dom}{{\rm dom}}
\newcommand{\cl}{{\rm cl}}
\newcommand{\id}{{\rm id}}

\renewcommand{\iff}{\mbox{\ \ \underline{\rm iff}\ \ }}

\newcommand{\forces}{\Vdash} 
\newcommand{\lesdot}{\mathrel{\mathord{<}\!\!\raise 0.8
pt\hbox{$\scriptstyle\circ$}}}  

\newcommand{\QED}{\hfill\vrule width 6pt height 6pt depth 0pt 
\vspace{0.1in}} 
\newcommand{\Proof}{\noindent {\sc Proof} \hspace{0.2in}} 

\newcommand{\cH}{{\cal H}}
\newcommand{\hchi}{\big(\cH(\chi^*),\mathord{\in},\mathord{<^*_{\chi^*}}\big)} 

\newcommand{\pre}[2]{{}^{\textstyle #1}{#2}}
\newcommand{\conc}{^\frown\!}

\newcommand{\bk}{{\bf k}}
\newcommand{\bl}{{\bf l}}
\newcommand{\bm}{{\bf m}}
\newcommand{\KL}{{\cal K}{\cal L}}
\newcommand{\cchi}{{\Bbb C}_\chi}
\newcommand{\Prs}{{\rm Pr}^\bigstar}
\newtheorem{theorem}{Theorem}[section] 
\newtheorem{claim}{Claim}[theorem]
 
\newtheorem{proposition}[theorem]{Proposition}

\newtheorem{kalpro}[theorem]{Kalikow Problem}
\newtheorem{definition}[theorem]{Definition}

\newtheorem{conclusion}[theorem]{Conclusion}

\newtheorem{remark}[theorem]{Remark} 
\newtheorem{finrem}[theorem]{Final Remarks}

\makeatletter
\newcounter{abc}

\newcommand{\labelabc}{(\alph{abc})}
\newenvironment{abc}{\list{\labelabc}{\usecounter{abc}%
      \advance\@enumdepth \@ne 
	\def\makelabel##1{\hss\llap{##1}}}}{\endlist}
\makeatother

\title{On a problem of Steve Kalikow}
\author{{\bf Saharon Shelah}\thanks{\ \ Research supported by  ``The
Israel Science Foundation'' administered by The Israel Academy of Sciences and
Humanities.  Publication no 590}\\ 
Institute of Mathematics, \\
The Hebrew University of Jerusalem,\\
Jerusalem 91904, Israel\\
and\\
Department of Mathematics,\\
Rutgers University,\\
New Brunswick, NJ 08854, USA}
\date{\today}

\begin{document}
\baselineskip13.14 truept

\maketitle

\begin{abstract}
The Kalikow problem for a pair $(\lambda,\kappa)$ of cardinal numbers,
$\lambda >\kappa$ (in particular $\kappa=2$) is whether we can map the family
of $\omega$--sequences from $\lambda$ to the family of $\omega$--sequences
from $\kappa$ in a very continuous manner. Namely, we demand that for $\eta,
\nu\in\lom$ we have:\quad $\eta,\nu$ are almost equal if and only if their
images are.

We show consistency of the negative answer e.g.~for $\aleph_\omega$ but we
prove it for smaller cardinals. We indicate a close connection with the free
subset property and its variants. 
\end{abstract}

\nocite{Ka90}
\nocite{Mi91}
\nocite{Sh:b}
\nocite{Sh:g}
\nocite{Sh:76}
\nocite{Sh:110}
\nocite{Sh:124}
\nocite{Sh:481}

\section{Introduction}
In the present paper we are interested in the following property of
pairs of cardinal numbers:

\begin{definition}
Let $\lambda,\kappa$ be cardinals. We say that the pair $(\lambda,\kappa)$ has
the Kalikow property  (and then we write $\KL(\lambda,\kappa)$) if
\begin{quotation}
\noindent there is a sequence $\langle F_n: n<\omega\rangle$ of functions such
that 
\[F_n:\pre n \lambda\longrightarrow \kappa\qquad\mbox{ (for $n<\omega$)}\]
and if $F:\lom\longrightarrow\kom$ is given by
\[(\forall\eta\in\lom)(\forall n\in\omega)\big(F(\eta)(n)=F_n(\eta\rest
n)\big)\] 
{\em then} for every $\eta,\nu\in\lom$
\[(\forall^\infty n)(\eta(n)=\nu(n))\iff (\forall^\infty n)(F(\eta)(n)=
F(\nu)(n)).\]
\end{quotation}
\end{definition}
In particular we answer the following question of Kalikow:

\begin{kalpro}
\label{kalpro}
Is $\KL(\con,2)$ provable in ZFC?
\end{kalpro}

The Kalikow property of pairs of cardinals was studied in \cite{Ka90}. Several
results are known already. Let us mention some of them. First, one can easily
notice that
\[\KL(\lambda,\kappa)\ \&\ \lambda'\leq\lambda\ \&\ \kappa'\geq\kappa\ \ \
\Rightarrow\ \ \ \KL(\lambda',\kappa').\]
Also (``transitivity'')
\[\KL(\lambda_2,\lambda_1)\ \&\ \KL(\lambda_1,\lambda_0)\ \ \ \Rightarrow\ \ \
\KL(\lambda_2,\lambda_0)\]
and
\[\KL(\lambda,\kappa)\ \ \ \Rightarrow\ \ \ \lambda\leq\kappa^{\aleph_0}.\]
Kalikow proved that CH implies $\KL(\con,2)$ (in fact that $\KL(\aleph_1,2)$
holds true) and he conjectured that CH is equivalent to $\KL(\con,2)$. 

The question \ref{kalpro} is formulated in \cite{Mi91} (Problem 15.15, p.
653).  

We shall prove that $\KL(\lambda,2)$ is closely tied with some variants of the
free subset property (both positively and negatively). First we present an
answer to the problem \ref{kalpro} proving the consistency of $\neg\KL(\con,
2)$ in \ref{main} (see \ref{finrem} too). Later we discuss variants of the
proof (concerning the cardinal and the forcing). Then we deal with positive
answer, in particular $\KL(\aleph_n,2)$ and we show that the negation of a
relative of the free subset property for $\lambda$ implies $\KL(\lambda,2)$.

We thank the participants of the Jerusalem Logic Seminar 1994/95 and
particularly Andrzej Ros{\l}anowski for writing it up so nicely.
\medskip

\noindent{\bf Notation:}\hspace{0.15in} We will use Greek letters
$\kappa,\lambda,\chi$ to denote (infinite) cardinals and letters
$\alpha,\beta,\gamma,\zeta,\xi$ to denote ordinals.  Sequences of ordinals
will be called $\bar{\alpha}$, $\bar{\beta}$, $\bar{\zeta}$ with the usual
convention that $\bar{\alpha}=\langle \alpha_n: n<\lh(\bar{\alpha})\rangle$
etc. Sets of ordinals will be denoted by $u$, $v$, $w$ (with possible
indexes). \\
The quantifiers $(\forall^\infty n)$ and $(\exists^\infty n)$ are
abbreviations for ``for all but finitely many $n\in\omega$'' and ``for
infinitely many $n\in\omega$'', respectively.

\section{The negative result}
For a cardinal $\chi$, the forcing notion $\cchi$ for adding $\chi$ many
Cohen reals consists of finite functions $p$ such that for some $w\in
[\chi]^{\textstyle {<}\omega}$, $n<\omega$ 
\[\dom(p)=\{(\zeta,k): \zeta\in w\ \&\ k<n\}\quad\mbox{ and }\quad\rng(p)
\subseteq 2\]
ordered by the inclusion.

\begin{theorem}
\label{main}
Assume $\lambda\rightarrow (\omega_1\cdot\omega)^{<\omega}_{2^\kappa}$,
$2^\kappa<\lambda\leq\chi$. Then
\[\forces_{\cchi} \neg\KL(\lambda,\kappa)\quad\mbox{ and hence }\quad
\forces_{\cchi}\neg\KL(\con,2).\] 
\end{theorem}

\Proof Suppose that $\cchi$-names $\name F_n$ (for $n\in\omega$) and a
condition $p\in\cchi$ are such that
\[p\forces_{\cchi}\mbox{``} \langle \name F_n: n<\omega\rangle\mbox{
examplifies } \KL(\lambda,\kappa)\mbox{''}.\]
For $\bar{\alpha}\in\pre n \lambda$ choose a maximal antichain $\langle
p_{\bar{\alpha},\ell}^n: \ell<\omega\rangle$ of $\cchi$ deciding the values of
$\name F_n(\bar{\alpha})$. Thus we have a sequence $\langle
\gamma^n_{\bar{\alpha},\ell}: \ell<\omega\rangle\subseteq\kappa$ such that
\[p^n_{\bar{\alpha},\ell}\forces_{\cchi}\name F_n(\bar{\alpha})=
\gamma^n_{\bar{\alpha},\ell}.\]
Let $\chi^*$ be a sufficiently large regular cardinal. Take an elementary
submodel $M$ of $\hchi$ such that
\begin{quotation}
\noindent $\|M\|=\chi$, $\chi+1\subseteq M$,

\noindent $\langle p_{\bar{\alpha},\ell}^n: \ell<\omega, n\in\omega,
\bar{\alpha}\in\pre n \lambda\rangle, \langle\gamma^n_{\bar{\alpha},\ell}:
\ell<\omega, n\in\omega, \bar{\alpha}\in\pre n \lambda\rangle\in M$.
\end{quotation}
By $\lambda\rightarrow (\omega_1\cdot\omega)^{<\omega}_{2^\kappa}$ (see
\cite{Sh:481}, Claim 1.3), we find a set $B\subseteq\lambda$ of indescernibles
in $M$ over  
\[\kappa\cup \{\langle
p^n_{\bar{\alpha},\ell}: \ell<\omega: n\in\omega,\bar{\alpha}\in\pre n
\lambda\rangle,\langle \gamma^n_{\bar{\alpha},\ell}: \ell<\omega:
n\in\omega,\bar{\alpha}\in\pre n \lambda\rangle, \chi,p\}\] 
and a system $\langle N_u: u\in [B]^{\textstyle {<}\omega}\rangle$ of
elementary submodels of $M$ such that 
\begin{abc}
\item $B$ is of the order type $\omega_1\cdot\omega$ and for $u,v\in
[B]^{\textstyle {<}\omega}$:  
\item $\kappa+1\subseteq N_u$, 
\item $\chi,p,\langle p^n_{\bar{\alpha},\ell}: \ell<\omega, n<\omega,
\bar{\alpha}\in\pre n \lambda\rangle, \langle\gamma^n_{\bar{\alpha},\ell}:
\ell<\omega, n<\omega, \bar{\alpha}\in\pre n \lambda\rangle\in N_u$, 
\item $|N_u|=\kappa$, $N_u\cap B = u$,
\item $N_u\cap N_v =N_{u\cap v}$,
\item $|u|=|v|\ \ \ \Rightarrow\ \ \ N_u\cong N_v$, and let $\pi_{u,v}:N_v
\longrightarrow N_u$ be this (unique) isomorphism,
\item $\pi_{v,v}=\id_{N_v}$, $\pi_{u,v}(v)=u$, $\pi_{u_0,u_1}\circ\pi_{u_1,
u_2}=\pi_{u_0,u_2}$, 
\item if $v'\subseteq v$, $|v|=|u|$ and $u'=\pi_{u,v}(v')$ then $\pi_{u',v'}
\subseteq\pi_{u,v}$.
\end{abc}
Note that if $u\subseteq B$ is of the order type $\omega$ then we may define
\[N_u=\bigcup\{N_v: v\mbox{ is a finite initial segment of }u\}.\]
Then the models $N_u$ (for $u\subseteq B$ of the order type $\leq\omega$) have
the properties (b)--(h) too.

Let $\langle\beta_\zeta: \zeta<\omega_1\cdot\omega\rangle$ be the increasing
enumeration of $B$. For a set $u\subseteq B$ of the order type $\leq\omega$
let $\bar{\beta}^u$ be the increasing enumeration of $u$ (so
$\lh(\bar{\beta}^u)=|u|$). Let $u^*=\{\beta_{\omega_1\cdot n}: n<\omega\}$.
For $k\leq\omega$ and a sequence $\bar{\xi}=\langle\xi_m:
m<k\rangle\subseteq\omega_1$ we define  
\[u[\bar{\xi}]=\{\beta_{\omega_1\cdot m + \xi_m}: m<k\}\cup
\{\beta_{\omega_1\cdot n}: n\in\omega\setminus k\}.\]
Now, working in $\V^{\cchi}$, we say that a sequence $\name{\bar{\xi}}$ is
{\em $k$--strange} if  
\begin{enumerate}
\item $\name{\bar{\xi}}$ is a sequence of countable ordinals greater than
$0$, $\lh(\name{\bar{\xi}})=k$
\item $(\forall m<\omega)(\name F_{m}(\bar{\beta}^{u[\name{\bar{\xi}}]}\rest
m)=\name F_{m}(\bar{\beta}^{u^*}\rest m))$.   
\end{enumerate}

\begin{claim}
\label{cl1}
In $\V^{\cchi}$:\\
{\em if} $\name{\bar{\xi}}^k$ are $k$--strange sequences (for $k<\omega$)
such that $(\forall k<\omega)(\name{\bar{\xi}}^k\vartriangleleft
\name{\bar{\xi}}^{k+1})$\\
{\em then} the sequence $\name{\bar{\xi}}\stackrel{\rm def}{=}
\bigcup\limits_{k<\omega}\name{\bar{\xi}}^k$ is $\omega$--strange.
\end{claim}

\noindent{\em Proof of the claim:}\hspace{0.15in} Should be clear (note that
in this situation we have $\bar{\beta}^{u[\name{\bar{\xi}}]}\rest
m=\bar{\beta}^{u[\name{\bar{\xi}}^m]}\rest m$). $\QED_{\ref{cl1}}$

\begin{claim}
\label{cl2}
\[p\forces_{\cchi}\mbox{``there are no $\omega$--strange sequences''}.\]
\end{claim}

\noindent{\em Proof of the claim:}\hspace{0.15in} Assume not. Then we find
a name $\name{\bar{\xi}}= \langle \name \xi_m: m<\omega\rangle$ for
an $\omega$--sequence and a condition $q\geq p$ such that
\[q\forces_{\cchi}\mbox{``}(\forall m<\omega)(0<\name\xi_m<\omega_1\quad
\&\quad\name F_m(\bar{\beta}^{u[\name{\bar{\xi}}]}\rest m) 
=\name F_m(\bar{\beta}^{u^*}\rest m))\mbox{''.}\]
By the choice of $p$ and $\name F_m$ we conclude that 
\[q\forces_{\cchi}\mbox{``}(\forall^\infty m)(\bar{\beta}^{u[
\name{\bar{\xi}}]}(m)=\bar{\beta}^{u^*}(m))\mbox{''}\] 
which contradicts the definition of $\bar{\beta}^{u[\name{\bar{\xi}}]}$,
$\bar{\beta}^{u^*}$ and the fact that 
\[q\forces_{\cchi}\mbox{``}(\forall m<\omega)(0<\name\xi_m<\omega_1)
\mbox{''}.\qquad\QED_{\ref{cl2}}\] 
\medskip

By \ref{cl1}, \ref{cl2}, any inductive attempt to construct (in $\V^{\cchi}$)
an $\omega$--strange sequence $\name{\bar{\xi}}$ has to fail. Consequently we
find a condition $p^*\geq p$, an integer $k<\omega$ and a sequence $\bar{\xi}=
\langle\xi_\ell: \ell<k\rangle$ such that
\[p^*\forces_{\cchi}\mbox{``}\bar{\xi}\mbox{ is $k$--strange but }
\neg(\exists\xi<\omega_1)(\bar{\xi}\conc\langle \xi\rangle\mbox{ is
$(k+1)$--strange})\mbox{''}.\]
Then in particular
\begin{description}
\item[$(\boxtimes)$]\hspace{0.15in} $p^*\forces_{\cchi}$ ``$(\forall
m<\omega)(\name F_m(\bar{\beta}^{u[\bar{\xi}]}\rest m)=\name
F_m(\bar{\beta}^{u^*}\rest m))$''. 
\end{description}
[It may happen that $k=0$, i.e. $\bar{\xi}=\langle\rangle$.]

\noindent For $\xi<\omega_1$ let $u_\xi=u[\bar{\xi}\conc\langle\xi\rangle]$ and
$w_\xi=u_\xi\cup (u^*\setminus\{\omega_1\cdot k\})$. Thus $w_0=u[\bar{\xi}]
\cup u^*$ and all $w_\xi$ have order type $\omega$ and
$\pi_{w_{\xi_1},w_{\xi_2}}$ is the identity on $N_{w_\xi\setminus\{\omega_1
\cdot k+\xi_2\}}$.\\ 
Let $q\stackrel{\rm def}{=}p^*\rest N_{w_0}$ and $q_\xi=\pi_{w_\xi,w_0}(q)\in
N_{w_\xi}$ (so $q_0=q$). As the isomorphism $\pi_{w_\xi,w_0}$ is the identity
on $N_{w_0}\cap N_{w_\xi}=N_{w_0\cap w_\xi}$ (and by the definition of Cohen
forcing), we have that the conditions $q,q_\xi$ are compatible. Moreover, as
$p^*\geq p$ and $p\in N_{\emptyset}$, we have that both $q$ and $q_\xi$ are
stronger than $p$. 

Now fix $\xi_0\in (0,\omega_1)$ (e.g. $\xi_0=1$) and look at the
sequences $\bar{\beta}^{u_{\xi_0}}$ and $\bar{\beta}^{u^*}$. They are
eventually equal and hence 
\[p\forces_{\cchi}\mbox{``}(\forall^\infty m)(\name
F_m(\bar{\beta}^{u_{\xi_0}}\rest m)=\name F_m(\bar{\beta}^{u^*} \rest
m))\mbox{''}.\] 
So we find $m^*<\omega$ and a condition $q^\prime_{\xi_0}\geq q_{\xi_0},
q$ such that
\begin{description}
\item[$(\otimes^{\xi_0,m^*}_{q^\prime_{\xi_0}})$] $q^\prime_{\xi_0}
\forces_{\cchi}$``$(\forall m\geq m^*)(\name F_m(\bar{\beta}^{u_{\xi_0}}
\rest m)=\name F_m(\bar{\beta}^{u^*}\rest m))$''

and (as we can increase $q^\prime_{\xi_0}$)
\item[$(\oplus^{\xi_0,m^*}_{q^\prime_{\xi_0}})$] the condition
$q^\prime_{\xi_0}$ decides the values of $\name
F_m(\bar{\beta}^{u_{\xi_0}}\rest m)$ and $\name F_m(\bar{\beta}^{u^*}\rest 
m)$ for all $m\leq m^*$.
\end{description}
Note that the condition $(\otimes^{\xi_0,m^*}_{q^\prime_{\xi_0}})$ means that
\begin{quotation}
\noindent there are NO $m\geq m^*$, $\ell_0,\ell_1<\omega$ with

\noindent $\gamma^m_{\bar{\beta}^{u_{\xi_0}}\rest m, \ell_0}\neq
\gamma^m_{\bar{\beta}^{u^*}\rest m, \ell_1}$ and the three conditions
$q^\prime_{\xi_0}$, $p^m_{\bar{\beta}^{u_{\xi_0}}\rest m, \ell_0}$ and
$p^m_{\bar{\beta}^{u^*}\rest m, \ell_1}$ have a common upper bound in $\cchi$
\end{quotation}
(remember the choice of the $p^n_{\bar{\alpha},\ell}$'s and
$\gamma^n_{\bar{\alpha},\ell}$'s). Similarly, the condition $(\oplus^{\xi_0,
m^*}_{q^\prime_{\xi_0}})$ means
\begin{quotation}
\noindent there are NO $m\leq m^*$, $\ell_0,\ell_1<\omega$ with

\noindent {\em either} $\gamma^m_{\bar{\beta}^{u_{\xi_0}}\rest m, \ell_0}\neq
\gamma^m_{\bar{\beta}^{u_{\xi_0}}\rest m, \ell_1}$ and both $q^\prime_{\xi_0}$
and $p^m_{\bar{\beta}^{u_{\xi_0}}\rest m,\ell_0}$, and $q^\prime_{\xi_0}$ and
$p^m_{\bar{\beta}^{u_{\xi_0}}\rest m,\ell_1}$ are compatible in $\cchi$ 

\noindent {\em or} $\gamma^m_{\bar{\beta}^{u^*}\rest m, \ell_0}\neq
\gamma^m_{\bar{\beta}^{u^*}\rest m, \ell_1}$ and both $q^\prime_{\xi_0}$ and
$p^m_{\bar{\beta}^{u^*}\rest m,\ell_0}$, and $q^\prime_{\xi_0}$ and
$p^m_{\bar{\beta}^{u^*}\rest m,\ell_1}$ are compatible in $\cchi$.
\end{quotation}
Consequently the condition $q^*_{\xi_0}\stackrel{\rm def}{=}
q^\prime_{\xi_0}\rest N_{w_0\cup w_{\xi_0}}$ has both properties 
$(\otimes^{\xi_0,m^*}_{q^*_{\xi_0}})$ and $(\oplus^{\xi_0,m^*}_{q^*_{\xi_0}})$
(and it is stronger than both $q$ and $q_{\xi_0}$).

Now, for $0<\xi<\omega_1$ let 
\[q^*_\xi\stackrel{\rm def}{=} \pi_{w_0\cup w_\xi, w_0\cup w_{\xi_0}}
(q^*_{\xi_0})\in N_{w_0\cup w_\xi}.\]
Then (for $\xi\in (0,\omega_1)$) the condition $q^*_\xi$ is
stronger than 
\[\mbox{both }\ q=\pi_{w_0\cup w_\xi, w_0\cup w_{\xi_0}}(q)\ \mbox{ and }\ 
q_\xi=\pi_{w_0\cup w_\xi, w_0\cup w_{\xi_0}}(q_{\xi_0})\]
and it has the properties $(\otimes^{\xi,m^*}_{q^*_{\xi}})$ and
$(\oplus^{\xi,m^*}_{q^*_{\xi}})$. Moreover for all $\xi_1,\xi_2$ the
conditions $q^*_{\xi_1}, q^*_{\xi_2}$ are compatible. [Why? By the definition
of Cohen forcing, and $\pi_{w_0\cup w_{\xi_2},w_0\cup w_{\xi_1}}(q_{\xi_1}^*)=
q_{\xi_2}^*$ (chasing arrows) and $\pi_{w_0\cup w_{\xi_2}, w_0\cup w_{\xi_1}}$
is the identity on $N_{w_0\cup w_{\xi_2}}\cap N_{w_0\cup w_{\xi_1}}=
N_{(w_0\cup w_{\xi_2})\cap (w_0\cup w_{\xi_1})}$ (see clauses (e), (f), (h)
above).]  

\begin{claim}
\label{cl3}
For each $\xi_1,\xi_2\in (0,\omega_1)$ the condition $q^*_{\xi_1}\cup
q^*_{\xi_2}$ forces in ${\cchi}$ that
\[(\forall m<\omega)(\name F_m(\bar{\beta}^{u_{\xi_1}}\rest m)= \name
F_m(\bar{\beta}^{u_{\xi_2}}\rest m)).\]
\end{claim}

\noindent{\em Proof of the claim:}\hspace{0.15in} If $m\geq m^*$ then, by
$(\otimes^{\xi_1,m^*}_{q^*_{\xi_1}})$, $(\otimes^{\xi_2,m^*}_{q^*_{\xi_2}})$
(passing through $\name F(\bar{\beta}^{u^*}\rest m)$) we get
\[q^*_{\xi_1}\cup q^*_{\xi_2}\forces_{\cchi}\mbox{``}\name
F_m(\bar{\beta}^{u_{\xi_1}}\rest m)= \name
F_m(\bar{\beta}^{u_{\xi_2}}\rest m)\mbox{''}.\]
If $m<m^*$ then we use $(\oplus^{\xi_1,m^*}_{q^*_{\xi_1}})$ and
$(\oplus^{\xi_1,m^*}_{q^*_{\xi_2}})$ and the isomorphism: the values assigned
by $q^*_{\xi_1}$, $q^*_{\xi_2}$ to $\name F_m(\bar{\beta}^{u_{\xi_1}}\rest m)$
and $\name F_m(\bar{\beta}^{u_{\xi_2}}\rest m)$ have to be equal (remember
$\kappa\subseteq N_\emptyset$, so the isomorphism is the identity on
$\kappa$). $\QED_{\ref{cl3}}$
\medskip

Look at the conditions
\[q_{\xi_1,\xi_2}\stackrel{\rm def}{=} q^*_{\xi_1}\rest
N_{w_{\xi_1}} \cup q^*_{\xi_2}\rest N_{w_{\xi_2}}\in N_{w_{\xi_1}
\cup w_{\xi_2}}.\]
It should be clear that for each $\xi_1,\xi_2\in (0,\omega_1)$ 
\[q_{\xi_1,\xi_2}\forces_{\cchi}\mbox{``}(\forall m<\omega)(\name
F_m(\bar{\beta}^{u_{\xi_1}}\rest m)= \name
F_m(\bar{\beta}^{u_{\xi_2}}\rest m))\mbox{''}.\]
Now choose $\xi\in (0,\omega_1)$ so large that
\[\dom(p^*)\cap (N_{w_\xi}\setminus N_{w_0}) =\emptyset\]
(possible as $\dom(p^*)$ is finite, use (e)). Take any $0<\xi_1<\xi_2<
\omega_1$ and put 
\[q^*\stackrel{\rm def}{=} \pi_{w_0\cup w_\xi,w_{\xi_1}\cup
w_{\xi_2}}(q_{\xi_1,\xi_2}).\]
(Note: $\pi_{w_0, w_{\xi_1}}\subseteq \pi_{w_0\cup w_\xi, w_{\xi_1}\cup
w_{\xi_2}}$ and $\pi_{w_\xi,w_{\xi_2}}\subseteq \pi_{w_0\cup w_\xi,
w_{\xi_1}\cup w_{\xi_2}}$.) By the isomorphism we get that 
\[q^*\forces_{\cchi}\mbox{``}(\forall m<\omega)(\name
F_m(\bar{\beta}^{u_{\xi}}\rest m)=\name F_m(\bar{\beta}^{u[\bar\xi]}\rest
m))\mbox{''}.\] 
Now look back:
\[\begin{array}{ll}
q^*_{\xi_1}\geq q_{\xi_1}=& \pi_{w_0\cup w_{\xi_1},w_0\cup w_{\xi_0}}
(q_{\xi_0})= \pi_{w_{\xi_1},w_{\xi_0}}(q_{\xi_0})=\\
\ &=\pi_{w_{\xi_1},w_{\xi_0}}(\pi_{w_{\xi_0},w_0}(q))=
\pi_{w_{\xi_1},w_0}(q)\\
\end{array}\]
and hence
\[q^*_{\xi_1}\rest N_{w_{\xi_1}}\geq \pi_{w_{\xi_1},w_0}(q)\]
and thus
\[q^*\rest N_{w_0}\geq\pi_{w_0,w_{\xi_1}}(q^*_{\xi_1}\rest N_{w_{\xi_1}})\geq
q=p^*\rest N_{w_0}.\]
Consequently, by the choice of $\xi$, the conditions $q^*$ and $p^*$ are
compatible (remember the definition of $q_{\xi_1,\xi_2}$ and $q^*$). Now use
$(\boxtimes)$ to conclude that 
\[q^*\cup p^*\forces_{\cchi}\mbox{``}(\forall m<\omega)(\name
F_m(\bar{\beta}^{u^*}\rest m)=\name F_m(\bar{\beta}^{u[\bar{\xi}]}\rest m)
=\name F_m(\bar{\beta}^{u_\xi}\rest m))\mbox{''}\]
which implies that
\[q^*\cup p^*\forces_{\cchi}\mbox{``}\bar{\xi}\conc\langle\xi\rangle\mbox{ is
$(k+1)$--strange''},\]
a contradiction. $\QED_{\ref{main}}$

\begin{remark}
About the proof of \ref{main}:
{\em
\begin{enumerate}
\item No harm is done by forgetting $0$ and replacing it by $\xi_1$, $\xi_2$.
\item A small modification of the proof shows that in $\V^{\cchi}$:
\begin{quotation}
{\em
\noindent If $F_n:\pre n\lambda\longrightarrow \kappa$ ($n\in\omega$) are such
that 
\[\hspace{-1cm}(\forall \eta,\nu\in\pre\omega\lambda)[(\forall^\infty
n)(\eta(n)=\nu(n))\ \ \Rightarrow\ \ (\forall^\infty n)(F_n(\eta\rest
n)=F_n(\nu\rest n))]\] 
then there are infinite sets $X_n\subseteq\lambda$ (for $n<\omega$) such that
\[(\forall n<\omega)(\forall \nu,\eta\in\prod_{\ell<n}X_\ell)(F_n(\nu)=
F_n(\eta)).\]
}
\end{quotation}
Say we shall have $X_n=\{\gamma_{n,i}: i<\omega\}$. Starting we have
$\gamma^*_0,\dots,\gamma^*_n,\ldots$. In the proof at stage $n$ we have
determined $\gamma_{\ell,i}$ ($\ell,i<n$) and $p\in G$, $p\in N_{\{\gamma_{
\ell,i}: \ell,i<\omega\}\cup\{\gamma^*_n,\gamma^*_{n+1},\ldots\}}$. For $n=0,1,
2$ as before. For $n+1>2$ first $\gamma_{0,n},\ldots,\gamma_{n-1,n}$ are easy
by transitivity of equalities. Then find $\gamma_{n,0},\gamma_{n,1}$ as before
then again duplicate. 
\item In the proof it is enough to use $\{\beta_{\omega\cdot n+\ell}:
n<\omega, \ell<\omega\}$. Hence, by 1.2 of \cite{Sh:481} it is enough to
assume $\lambda\rightarrow (\omega^3)^{<\omega}_{2^\kappa}$. This condition
is compatible with $\V={\bf L}$.
\item We can use only $\lambda\rightarrow (\omega^2)^{<\omega}_{2^\kappa}$. 
\end{enumerate}
}
\end{remark}

\begin{definition} 
\label{really}
\begin{enumerate}
\item For a sequence $\bar{\lambda}=\langle\lambda_n: n<\omega\rangle$ of
cardinals we define the property $(\circledast)_{\bar{\lambda}}$:
\begin{description}
\item[$(\circledast)_{\bar{\lambda},}$] for every model $M$ of a
countable language, with universe $\sup\limits_{n\in\omega}\lambda_n$ and
Skolem functions (for simplicity) there is a sequence $\langle X_n:
n<\omega\rangle$ such that
\begin{abc}
\item $X_n\in [\lambda_n]^{\textstyle \lambda_n}$ (actually $X_n\in
[\lambda_n]^{\textstyle \omega_1}$ suffices)
\item for every $n<\omega$ and $\bar{\alpha}=\langle \alpha_\ell: \ell\in
[n+1,\omega)\rangle\in \prod\limits_{\ell\geq n+1} X_\ell$, letting (for
$\xi\in X_{n}$) 
\[M^\xi_{\bar{\alpha}}={\rm Sk}(\bigcup_{\ell<n} X_\ell\cup\{\xi\}\cup
\{\alpha_\ell: \ell\in [n+1,\omega)\})\]
we have:
\begin{description}
\item[$(\bigoplus)$] the sequence $\langle M^\xi_{\bar{\alpha}}: \xi\in
X_n\rangle$ forms a $\Delta$--system with the heart $N_{\bar{\alpha}}$ and 
its elements are pairwise isomorphic over the heart $N_{\bar{\alpha}}$.
\end{description}
\end{abc}
\end{description}
\item For a cardinal $\lambda$ the condition $(\circledast)^\lambda$ is:
\begin{description}
\item[$(\circledast)^\lambda$] there exists a sequence $\bar{\lambda}=
\langle\lambda_n: n<\omega\rangle$ such that $\sum\limits_{n<\omega}\lambda_n
= \lambda$ and the condition $(\circledast)_{\bar{\lambda}}$ holds true.
\end{description}
\end{enumerate}
\end{definition}

In \cite{Sh:76} a condition $(*)_{\lambda}$, weaker than
$(\circledast)^\lambda$ was considered. Now, \cite{Sh:124} continues
\cite{Sh:76} to get stronger indiscernibility. But by the same proof (using
$\omega$-measurable) one can show the consistency of
$(\circledast)^{\aleph_\omega} + {\rm GCH}$.  

Now note that to carry the proof of \ref{main} we need even less then
$(\circledast)^\lambda$: the $\bigcup\limits_{\ell<n} X_\ell$ (in (b) of
\ref{really}) is much more then needed; it suffices to have
$\bar{\beta}^0\cup\bar{\beta}^1$ where $\bar{\beta}^0,\bar{\beta}^1\in
\prod\limits_{\ell<n} X_\ell$.

\begin{conclusion}
It is consistent that
\[\con=\aleph_{\omega+1}\quad\mbox{ and }\quad \bigwedge_{n<\omega}
\neg \KL(\aleph_\omega,\aleph_n) \quad\mbox{ so $\neg\KL(\con,2)$}.\]
\end{conclusion}

\begin{remark}
\label{koepke}
{\em
Koepke \cite{Ko84} continues \cite{Sh:76} to get equiconsistency. His
refinement of \cite{Sh:76} (for the upper bound) works below too. 
}
\end{remark}

\section{The positive result}
For an algebra $M$ on $\lambda$ and a set $X\subseteq\lambda$ the closure of
$X$ under functions of $M$ is denoted by $\cl_M(X)$. Before proving our result
(\ref{second}) we remind the reader of some definitions and propositions. 

\begin{proposition}
\label{stars}
For an algebra $M$ on $\lambda$ the following conditions are equivalent 
\begin{description}
\item[$(\bigstar)^0_M$] for each sequence $\langle\alpha_n: n\in\omega\rangle
\subseteq\lambda$ we have
\[(\forall^\infty n)(\alpha_n\in\cl_M(\{\alpha_k: n<k<\omega\})),\]
\item[$(\bigstar)^1_M$] there is no sequence $\langle A_n: n\in\omega\rangle
\subseteq [\lambda]^{\aleph_0}$ such that
\[(\forall n\in \omega)(\cl_M(A_{n+1})\varsubsetneq\cl_M(A_n)),\]
\item[$(\bigstar)^2_M$] $(\forall A\in [\lambda]^{\aleph_0})(\exists B\in
[A]^{\aleph_0})(\forall C\in [B]^{\aleph_0})(\cl_M(B)=\cl_M(C))$.
$\QED_{\ref{stars}}$ 
\end{description}
\end{proposition}

\begin{definition}
We say that a cardinal $\lambda$ has the $(\bigstar)$--property for $\kappa$
(and then we write $\Prs(\lambda,\kappa)$) if there is an algebra $M$ on
$\lambda$ with vocabulary of cardinality $\leq\kappa$ satisfying one
(equivalently: all) of the conditions $(\bigstar)^i_M$ ($i<3$) of \ref{stars}.
If $\kappa=\aleph_0$ we may omit it.
\end{definition}

Remember
\begin{proposition}
\label{V0V1}
If $\V_0\subseteq\V_1$ are universes of set theory, $\V_1\models
\neg\Prs(\lambda)$ then $\V_0\models\neg\Prs(\lambda)$. 
\end{proposition}

\Proof By absoluteness of the existence of an $\omega$--branch to a tree.
$\QED_{\ref{V0V1}}$ 

\begin{remark}
{\em
The property $\neg\Prs(\lambda)$ is a kind of a large cardinal property. It
was clarified in {\bf L} (remember that it is inherited from $\V$ to {\bf L})
by Silver \cite{Si70} to be equiconsistent with ``there is a beautiful
cardinal'' (terminology of 2.3 of \cite{Sh:110}), another partition property
inherited by {\bf L}. 
}
\end{remark}

\begin{proposition}
\label{getPn}
For each $n\in\omega$, $\Prs(\aleph_n)$.
\end{proposition}

\Proof This was done in chapter XIII of \cite{Sh:b}, see chapter VII of
\cite{Sh:g} too, and probably earlier by Silver. However, for the sake of
completeness we will give the proof. 

First note that clearly $\Prs(\aleph_0)$ and thus we have to deal with
the case when $n>0$. Let $f,g:\aleph_n\longrightarrow\aleph_n$ be two
functions such that 
\begin{quotation}
\noindent if $m<n$, $\alpha\in [\aleph_m,\aleph_{m+1})$ 

\noindent then $f(\alpha,\cdot)\rest\alpha: \alpha\stackrel{\rm
1-1}{\longrightarrow}\aleph_m$, $g(\alpha,\cdot)\rest\aleph_m:
\aleph_m\stackrel{\rm 1-1}{\longrightarrow}\alpha$ are functions inverse each
to the other. 
\end{quotation}
Let $M$ be the following algebra on $\aleph_n$:
\[M=(\aleph_n, f, g, m)_{m\in\omega}.\]
We want to check the condition $(\bigstar)^1_M$:\\
assume that a sequence $\langle A_k: k<\omega\rangle\subseteq
[\aleph_n]^{\aleph_0}$ is such that for each $k<\omega$
\[\cl_M(A_{k+1})\varsubsetneq\cl_M(A_k).\]
For each $m<n$, the sequence $\langle\sup(\cl_M(A_k)\cap\aleph_{m+1}):
k<\omega\rangle$ is non-increasing and therefore it is eventually
constant. Consequently we find $k^*$ such that
\[(\forall m<n)(\sup(\cl_M(A_{k^*+1})\cap\aleph_{m+1})=\sup(\cl_M(A_{k^*})
\cap\aleph_{m+1})).\]
By the choice of $\langle A_k: k<\omega\rangle$ we have
$\cl_M(A_{k^*+1})\varsubsetneq\cl_M(A_{k^*})$. Let
\[\alpha_0\stackrel{\rm def}{=}\min(\cl_M(A_{k^*})\setminus
\cl_M(A_{k^*+1})).\]
As the model $M$ contains individual constants $m$ (for $m\in\omega$) we know
that $\aleph_0\subseteq\cl_M(\emptyset)$ and hence $\aleph_0\leq\alpha_0$. Let
$m<n$ be such that $\aleph_m\leq\alpha_0<\aleph_{m+1}$. By the choice of $k^*$
we find $\beta\in\cl_M(A_{k^*+1})\cap\aleph_{m+1}$ such that
$\alpha_0\leq\beta$. Then necessarily $\alpha_0<\beta$. Look at
$f(\beta,\alpha_0)$: we know that $\alpha_0,\beta\in\cl_M(A_{k^*})$ and
therefore $f(\beta,\alpha_0)\in\cl_M(A_{k^*})\cap\aleph_m$ and
$f(\beta,\alpha_0)<\alpha_0$. The minimality of $\alpha_0$ implies that
$f(\beta,\alpha_0)\in\cl_M(A_{k^*+1})$ and hence
\[\alpha_0=g(\beta, f(\beta,\alpha_0))\in \cl_M(A_{k^*+1}),\]
a contradiction. $\QED_{\ref{getPn}}$
\medskip

\noindent{\bf Explanation:}\quad Better think of the proof from the end. Let
$\bar{\alpha}=\langle\alpha_n:n<\omega\rangle\in\lom$. So for some $n(*)$,
$n(*)\leq n<\omega\ \Rightarrow\ \alpha_n\in\cl_M(\alpha_\ell: \ell>n)$. So
for some $m_n>n$, $\{\alpha_{n(*)},\ldots,\alpha_{n-1}\}\subseteq\cl_M(
\alpha_n,\ldots,\alpha_{m-1})$ and 
\[(\forall\ell<n(*))(\alpha_\ell\in\cl_M(\alpha_\ell:\ell>n(*))\ \Rightarrow\
\alpha_\ell\in\cl_M(\alpha_\ell:\ell\in [n,m_n))).\]
Let $W^*=\{\ell<n(*):\alpha_\ell\in\cl_M(\alpha_n:n\geq n(*))$. It is natural
to aim at:
\begin{description}
\item[$(*)$] for $n$ large enough (say $n>m_{n(*)}$), $F_n(\langle
\alpha_\ell:\ell<n\rangle)$ depends just on $\{\alpha_\ell:\ell\in [n(*),n)
\mbox{ or }\ell\in w\}$ and $\langle F_m(\bar{\alpha}\rest m): m\geq n\rangle$
codes $\bar{\alpha}\rest (w\cup [n(*),\omega))$.
\end{description}
Of course, we are a given $n$ and we do not know how to compute the real
$n(*)$, but we can approximate. Then we look at a late enoug end segment where
we compute down.

\begin{theorem}
\label{second}
Assume that $\lambda\leq\con$ is such that $\Prs(\lambda)$ holds.\\
Then $\KL(\lambda,\omega)$ (and hence $\KL(\lambda,2)$).
\end{theorem}

\Proof We have to construct functions $F_n:\pre n\lambda\longrightarrow\omega$
witnessing $\KL(\lambda,\omega)$. For this we will introduce functions $\bk$
and $\bl$ such that for $\bar{\alpha}\in\pre n\lambda$ the value of
$\bk(\bar{\alpha})$ will say which initial segment of $\bar{\alpha}$ will be
irrelevant for $F_n(\bar{\alpha})$ and $\bl(\bar{\alpha})$ will be such that
(under certain circumstances) elements $\alpha_i$ (for $\bk(\bar{\alpha})\leq
i<\bl(\bar{\alpha})$) will be encoded by $\langle\alpha_j: j\in
[\bl(\bar{\alpha}),n)\rangle$. 

Fix a sequence $\langle\eta_\alpha: \alpha<\lambda\rangle\subseteq\can$
with no repetitions.

Let $M$ be an algebra on $\lambda$ such that $(\bigstar)^0_M$ holds true. We
may assume that there are no individual constants in $M$ (so
$\cl_M(\emptyset)=\emptyset$).\\ 
Let $\langle \tau^n_\ell(x_0,\ldots,x_{n-1}): \ell<\omega\rangle$ list all
$n$-place terms of the language of the algebra $M$ (and $\tau^1_0(x)$ is
$x$). For $\bar{\alpha}\in\pre{\omega{\geq}}\lambda$ (with $\alpha_j$ the
$j$-th element in $\bar{\alpha}$) let 
\[u(\bar{\alpha})=\{\ell<\lh(\bar{\alpha}):\alpha_\ell\notin
\cl_M\big(\bar{\alpha}\rest (\ell,\lh(\bar{\alpha}))\big)\}\cup \{0\}\]
and for $\ell\notin u(\bar{\alpha})$, $\ell<\lh(\bar{\alpha})$ let
\[\begin{array}{lcl}
f_\ell(\bar{\alpha})&=&\min\{j: \alpha_\ell\in\cl_M(\bar{\alpha}\rest
(\ell,j))\}\\
g_\ell(\bar{\alpha})&=&\min\{i: \alpha_\ell=
\tau^{f_\ell(\bar{\alpha})-\ell-1}_i(\bar{\alpha}\rest
(\ell,f_\ell(\bar{\alpha})))\}.\\ 
  \end{array}\]
For $\bar{\alpha}\in\pre n\lambda$ ($1<n<\omega$) put 
\[\begin{array}{lcl}
k_1(\bar{\alpha})&=&\min\big((u(\bar{\alpha}\rest (n-1))\setminus
u(\bar{\alpha}))\cup\{n-1\}\big)\\
k_0(\bar{\alpha})&=&\max\big(u(\bar{\alpha})\cap k_1(\bar{\alpha})\big).\\
  \end{array}\]
Note that if ($n>1$ and) $\bar{\alpha}\in\pre n\lambda$ then $n-1\in
u(\bar{\alpha})$ (as $\cl_M(\emptyset)=\emptyset$) and $k_1(\bar{\alpha})>0$
(as always $0\in u(\bar{\beta})$) and $k_0(\bar{\alpha})$ is well defined (as
$0\in u(\bar{\alpha})\cap k_1(\bar{\alpha})$) and $k_0(\bar{\alpha})<
k_1(\bar{\alpha})<n$. Moreover, for all $\ell\in (k_0(\bar{\alpha}),
k_1(\bar{\alpha}))$ we have $\alpha_\ell\notin u(\bar{\alpha}\rest (n-1))$ and
thus $\alpha_\ell\in\cl_M(\bar{\alpha}\rest (\ell,n-1))$. Now, for
$\bar{\alpha}\in\pre{\omega{>}}\lambda$, $\lh(\bar{\alpha})>1$ we define  
\[\begin{array}{lcl}
\bl(\bar{\alpha})&=&\max\{j\leq k_1(\bar{\alpha}): j{>} k_0(\bar{\alpha})\
\ \Rightarrow\ \ (\forall i\!\in\! (k_0(\bar{\alpha}),j))(g_i(\bar{\alpha})\leq
\lh(\bar{\alpha}))\}\\  
\bm(\bar{\alpha})&=&\max\{j\leq\bl(\bar{\alpha}): j{>} \max\{1,
k_0(\bar{\alpha})\}\ \ \Rightarrow\ \ k_0(\bar{\alpha}\rest j)=
k_0(\bar{\alpha})\}\\ 
\bk(\bar{\alpha})&=&\bl(\bar{\alpha}\rest \bm(\bar{\alpha}))\ \
\mbox{ (if $\bm(\bar{\alpha})\leq 1$ then put $\bk(\bar{\alpha})=-1$)}.\\
  \end{array}\]
Clearly $\bk(\bar{\alpha})<\bm(\bar{\alpha})\leq\bl(\bar{\alpha})\leq
k_1(\bar{\alpha})<\lh(\bar{\alpha})$.

\begin{claim}
\label{cl4}
For each $\bar{\alpha}\in\pre\omega\lambda$, the set $u(\bar{\alpha})$ is
finite and:
\begin{enumerate}
\item The sequence $\langle k_1(\bar{\alpha}\rest n): n<\omega\rangle$
diverges to $\infty$.
\item The sequence $\langle
k_0(\bar{\alpha}\rest n): n<\omega\ \&\ k_0(\bar{\alpha})\neq \max
u(\bar{\alpha})\rangle$, if infinite, diverges to $\infty$. There are
infinitely many $n<\omega$ with $k_0(\bar{\alpha}\rest n)=\max
u(\bar{\alpha})$.  
\item The sequence $\langle\bl(\bar{\alpha}\rest n): n<\omega\rangle$ diverges
to $\infty$.
\item The sequences $\langle\bm(\bar{\alpha}\rest n): n<\omega\rangle$ and
$\langle\bk(\bar{\alpha}\rest n): n<\omega\rangle$ diverge to $\infty$.
\end{enumerate}
\end{claim}

\noindent{\em Proof of the claim:}\hspace{0.15in} Let $\bar{\alpha}=\langle
\alpha_n: n<\omega\rangle\in\pre \omega\lambda$. By the property
$(\bigstar)^0_M$ we find $n^*<\omega$ such that $u(\bar{\alpha})\subseteq
n^*$. Fix $n_0>n^*$ and define
\[n_{1}=\max\{f_n(\bar{\alpha})+g_n(\bar{\alpha})+2: n\in (n_0+1)\setminus
u(\bar{\alpha})\}\] 
(so $n_{1}\geq f_{n_0}(\bar{\alpha})+2>n_0+3$ and for all $\ell\in
(n_0+1)\setminus u(\bar{\alpha})$ we have:\\
$\alpha_\ell\in\cl_M(\alpha_{\ell+1},\ldots,\alpha_{n_{1}-1})$ is witnessed by
$\tau^{f_\ell(\bar{\alpha})-\ell-1}_{g_\ell(\bar{\alpha})}(\alpha_{\ell+1},
\ldots,\alpha_{f_\ell(\bar{\alpha})-1})$ with $f_\ell(\bar{\alpha}),
g_\ell(\bar{\alpha})<n_1-1$).     
\medskip

\noindent 1)\ \ \ Note that $u(\bar{\alpha}\rest n)\cap
(n_0+1)=u(\bar{\alpha})$ for all $n\geq n_1-1$ and hence for $n\geq n_1$
\[u(\bar{\alpha}\rest n)\cap (n_0+1)= u(\bar{\alpha}\rest (n-1))\cap (n_0+1).\]
Consequently for all $n\geq n_1$ we have that $k_1(\bar{\alpha}\rest n)>n_0$.
As we could have chosen $n_0$ arbitrarily large we may conclude that
$\lim\limits_{n\to\infty} k_1(\bar{\alpha}\rest n)=\infty$.
\medskip

\noindent 2)\ \ \ Note that for all $n\geq n_1$
\[\mbox{either }\ k_0(\bar{\alpha}\rest n)=\max(u(\bar{\alpha}))\ \mbox{ or }\
k_0(\bar{\alpha}\rest n)>n_0.\]
Hence, by the arbirarity of $n_0$, we get the first part of 2).\\  
Let $\ell^*=\min (u(\bar{\alpha}\rest n_1)\setminus u(\bar{\alpha}))$ (note
that $n_1-1\in u(\bar{\alpha}\rest n_1)\setminus u(\bar{\alpha})$). Clearly
$\ell^*>n_0$ and $\alpha_{\ell^*}\notin u(\bar{\alpha})$. Consider
$n=f_{\ell^*}(\bar{\alpha})$ (so $\ell^*\leq n-2$, $n_1\leq n-1$). Then
$\ell^*\in u(\bar{\alpha}\rest (n-1))\setminus u(\bar{\alpha}\rest 
n)$. As 
\[\ell^*\cap u(\bar{\alpha}\rest n_1)= \ell^*\cap u(\bar{\alpha}\rest n-1)=
u(\bar{\alpha})\] 
(remember the choice of $\ell^*$) we conclude that 
\[\ell^*=k_1(\bar{\alpha}\rest n)\quad\mbox{ and }\quad k_0(\bar{\alpha}\rest
n) = \max u(\bar{\alpha}).\]
Now, since $n_0$ was arbitrarily large, we conclude that for infinitely many
$n$, $k_0(\bar{\alpha}\rest n)=\max u(\bar{\alpha})$.
\medskip

\noindent 3)\ \ \ Suppose that $n\geq n_1$. Then we know that 
$k_1(\bar{\alpha}\rest n)> n_0$ and either $k_0(\bar{\alpha}\rest n)= \max
u(\bar{\alpha})$ or $k_0(\bar{\alpha}\rest n)> n_0$ (see above). If the first
possibility takes place then, as $n\geq n_1$, we may use $j=n_0+1$ to witness
that $\bl(\bar{\alpha}\rest n)>n_0$ (remember the choice of $n_1$). If
$k_0(\bar{\alpha}\rest n)> n_0$ then clearly $\bl(\bar{\alpha}\rest n)> n_0$.
As $n_0$ could be arbitrarily large we are done. 
\medskip

\noindent 4)\ \ \ Suppose we are given $m_0<\omega$. Take $m_1>m_0$ such that
for all $n\geq m_1$ 
\[\mbox{either } k_0(\bar{\alpha}\rest n) =\max u(\bar{\alpha})\ \mbox{ or }\
k_0(\bar{\alpha}\rest n)>m_0\]
(possible by 2)) and then choose $m_2>m_1$ such that $k_0(\bar{\alpha}\rest
m_2)=\max u(\bar{\alpha})$ (by 2)). Due to 3) we find $m_3>m_2$ such that for
all $n\geq m_3$, $\bl(\bar{\alpha}\rest n)>m_2$.\\
Now suppose that $n\geq m_3$. If $k_0(\bar{\alpha}\rest n)= \max
u(\bar{\alpha})$ then, as $\bl(\bar{\alpha}\rest n)>m_2$, we get 
$\bm(\bar{\alpha}\rest n)\geq m_2>m_0$. Otherwise $k_0(\bar{\alpha}\rest n)>
m_0$ (as $n>m_1$) and hence $\bm(\bar{\alpha}\rest n)>m_0$. This shows that
$\lim\limits_{n\to\infty}\bm(\bar{\alpha}\rest n)=\infty$. Now, immediately by
the definition of $\bk$ and 3) above we conclude that
$\lim\limits_{n\to\infty}\bk(\bar{\alpha}\rest n)=\infty$.  $\QED_{\ref{cl4}}$

\begin{claim}
\label{cl5}
If $\bar{\alpha}^1,\bar{\alpha}^2\in\pre\omega\lambda$ are such that
$(\forall^\infty n)(\alpha_n^1=\alpha_n^2)$ then 
\[(\forall^\infty n)\bigg(\bl(\bar{\alpha}^1\rest n)=\bl(\bar{\alpha}^2\rest
n)\ \&\ \bm(\bar{\alpha}^1\rest n)=\bm(\bar{\alpha}^2\rest n)\ \&\
\bk(\bar{\alpha}^1\rest n)=\bk(\bar{\alpha}^2\rest n)\bigg).\] 
\end{claim}

\noindent{\em Proof of the claim:}\hspace{0.15in} Let $n_0$ be greater than
$\max(u(\bar{\alpha}^1)\cup u(\bar{\alpha}^2))$ and such that 
\[\bar{\alpha}^1\rest [n_0,\omega) = \bar{\alpha}^2\rest [n_0,\omega).\]
For $k=1,2,3$ define $n_k$ by
\[n_{k+1}=\max\{f_n(\bar{\alpha}^i)+g_n(\bar{\alpha}^i)+2: n\in
(n_k+1)\setminus u(\bar{\alpha}^i),\ i<2\}.\]
As in the proof of \ref{cl4} we have that then for $i=1,2$ and $j<3$:
\begin{description}
\item[$(\otimes^1)$] $(\forall n\geq n_{j+1})(k_0(\bar{\alpha}^i\rest n)= \max
u(\bar{\alpha}^i)\quad\mbox{ or }\quad k_0(\bar{\alpha}^i\rest n)>n_j)$
\item[$(\otimes^2)$] $(\forall n\geq n_{j+1})(k_1(\bar{\alpha}^i\rest n)> n_j\
\&\ \bl(\bar{\alpha}^i\rest n)> n_j)$
\item[$(\otimes^3)$] $(\exists n'\in (n_1,n_2))(k_0(\bar{\alpha}^1\rest n')=
\max u(\bar{\alpha}^1)\ \&\ k_0(\bar{\alpha}^2\rest n')=\max
u(\bar{\alpha}^2))$
\end{description}
(for $(\otimes^3)$ repeat arguments from \ref{cl4}.(2) and use the fact that
$\bar{\alpha}^1\rest [n_0,\omega)=\bar{\alpha}^2\rest [n_0,\omega)$). Clearly
\begin{description}
\item[$(\otimes^4)$] $(\forall n>n_0)(u(\bar{\alpha}^1\rest n)\setminus n_0=
u(\bar{\alpha}^2\rest n) \setminus n_0)$.
\end{description}
Hence, applying $(\otimes^1)$, $(\otimes^2)$, we conclude that:
\begin{description}
\item[$(\otimes^5)$] $(\forall n\geq n_{1})(k_1(\bar{\alpha}^1\rest n) =
k_1(\bar{\alpha}^2\rest n))$ and
\item[$(\otimes^6)$] for all $n\geq n_1$:
\begin{quotation}
\noindent either $k_0(\bar{\alpha}^1\rest n)=\max u(\bar{\alpha}^1)$ and
$k_0(\bar{\alpha}^2\rest n)=\max u(\bar{\alpha}^2)$ 

\noindent or $k_0(\bar{\alpha}^1\rest n)=k_0(\bar{\alpha}^2\rest n)$.
\end{quotation}
\end{description}
Since
\[(\forall n\geq n_0)(f_n(\bar{\alpha}^1)=f_n(\bar{\alpha}^2)\ \&\
g_n(\bar{\alpha}^1)=g_n(\bar{\alpha}^2))\]
and by $(\otimes^2)+(\otimes^5)$, we get (compare the proof of \ref{cl4}):
\[(\forall n\geq n_1)(\bl(\bar{\alpha}^1\rest n)=\bl(\bar{\alpha}^2\rest n))\]
and by $(\otimes^2)+(\otimes^3)+(\otimes^6)$
\[(\forall n\geq n_3)(\bm(\bar{\alpha}^1\rest n)= \bm(\bar{\alpha}^2\rest
n)\geq n_1).\]
Moreover, now we easily get that 
\[(\forall n\geq n_3)(\bk(\bar{\alpha}^1\rest n)= \bk(\bar{\alpha}^2\rest
n)). \qquad\QED_{\ref{cl5}}\]
\medskip

\noindent For integers $n_0\leq n_1\leq n_2$ we define functions
$F^0_{n_0,n_1,n_2}: \pre {n_2} \lambda\longrightarrow \cH(\aleph_0)$ by
letting $F^0_{n_0,n_1,n_2}(\alpha_0,\ldots,\alpha_{n_2-1})$ (for
$\langle\alpha_0,\ldots,\alpha_{n_2-1}\rangle\in \pre{n_2}\lambda$) be the
sequence consisting of:
\begin{abc}
\item $\langle n_0,n_1,n_2\rangle$,
\item the set $T_{n_1,n_2}$ of all terms $\tau^n_\ell$ such that $n\leq
n_2-n_1$ and 
\begin{quotation}
\noindent either $\ell\leq n_2$ (we will call it {\em the simple case})

\noindent or $\tau^n_\ell$ is a composition of depth at most $n_2$ of such
terms,
\end{quotation}
\item $\langle \eta_{\alpha}\rest n_2, n,\ell,\langle i_0,\ldots,i_{n-1}
\rangle\rangle$ for $n\leq n_2-n_1$, $i_0,\ldots,i_{n-1}\in [n_1,n_2)$ and
$\ell$ such that $\tau^n_\ell\in T_{n_1,n_2}$ and $\alpha=\tau^n_\ell(
\alpha_{i_0}, \ldots, \alpha_{i_{n-1}})$,  
\item $\langle n,\ell,\langle i_0,\ldots,i_{n-1}\rangle, i\rangle$ for $n\leq
n_2-n_1$, $i_0,\ldots,i_{n-1}\in [n_1,n_2)$, $i\in [n_0,n_1)$ and $\ell$ such
that $\tau^n_\ell\in T_{n_1,n_2}$ and $\alpha_i=\tau^n_\ell(\alpha_{i_0},
\ldots, \alpha_{i_{n-1}})$, 
\item equalities among appropriate terms, i.e. all tuples
\[\langle n',\ell',n'',\ell'',\langle i'_0,\ldots,i'_{n'-1}\rangle, \langle
i''_0,\ldots,i''_{n''-1}\rangle\rangle\]
such that $n_1\leq i'_0<\ldots<i'_{n'-1}<n_2$, $n_1\leq
i''_0<\ldots<i''_{n''-1}<n_2$, $n',n''\leq n_2-n_1$, $\ell',\ell''$ are such
that $\tau^{n'}_{\ell'},\tau^{n''}_{\ell''}\in T_{n_1,n_2}$ and 
\[\tau^{n'}_{\ell'}(\alpha_{i'_0},\ldots,\alpha_{i'_{n'-1}})=
\tau^{n''}_{\ell''}(\alpha_{i''_0},\ldots,\alpha_{i''_{n''-1}}).\]
\end{abc}

\noindent (Note that the value of $F^0_{n_0,n_1,n_2}(\bar{\alpha})$ does not
depend on $\bar{\alpha}\rest n_0$.)

\noindent Finally we define functions $F_n:\pre n\lambda\longrightarrow
\cH(\aleph_0)$ (for $1<n<\omega$) by: 
\begin{quotation}
\noindent if $\bar{\alpha}\in\pre n\lambda$

\noindent then $F_n(\bar{\alpha})=F^0_{\bk(\bar{\alpha}),\bl(\bar{\alpha}),n}
(\bar{\alpha})$. 
\end{quotation}
As $\cH(\aleph_0)$ is countable we may think that these functions are into
$\omega$. We are going to show that they witness
$\KL(\lambda,\omega)$. 

\begin{claim}
\label{cl7}
If $\bar{\alpha}^1,\bar{\alpha}^2\in\pre\omega\lambda$ are such that
$(\forall^\infty n)(\alpha_n^1=\alpha^2_n)$\\
then $(\forall^\infty n)(F_n(\bar{\alpha}^1\rest n)=F_n(\bar{\alpha}^2\rest
n))$. 
\end{claim}

\noindent{\em Proof of the claim:}\hspace{0.15in} Take $m_0<\omega$ such that
for all $n\in [m_0,\omega)$ we have
\[\alpha^1_n=\alpha^2_n,\quad \bl(\bar{\alpha}^1\rest n)=\bl(\bar{\alpha}^2
\rest n),\quad\mbox{ and }\bk(\bar{\alpha}^1\rest n)=\bk(\bar{\alpha}^2\rest
n)\]
(possible by \ref{cl5}). Let $m_1> m_0$ be such that for all $n\geq m_1$:
\[\bk(\bar{\alpha}^1\rest n)=\bk(\bar{\alpha}^2\rest n) > m_0\]
(use \ref{cl4}). Then, for $n\geq m_1$, $i=1,2$ we have
\[F_n(\bar{\alpha}^i\rest n)= F^0_{\bk(\bar{\alpha}^i\rest n),
\bl(\bar{\alpha}^i\rest n), n}(\bar{\alpha}^i\rest
n)=F^0_{\bk(\bar{\alpha}^1\rest n),\bl(\bar{\alpha}^1\rest
n),n}(\bar{\alpha}^i\rest n).\]
Since the value of $F^0_{n_0,n_1,n_2}(\bar{\beta})$ does not depend on
$\bar{\beta}\rest n_0$ and the sequences $\bar{\alpha}^1\rest n$,
$\bar{\alpha}^2\rest n$ agree on $[m_0,\omega)$, we get 
\[F^0_{\bk(\bar{\alpha}^1\rest n),\bl(\bar{\alpha}^1\rest n), n}(\bar{\alpha}^1
\rest n)=F^0_{\bk(\bar{\alpha}^1\rest n), \bl(\bar{\alpha}^1\rest n), n}(
\bar{\alpha}^2\rest n)=F^0_{\bk(\bar{\alpha}^2\rest n),\bl(\bar{\alpha}^2\rest
n), n}(\bar{\alpha}^2 \rest n),\]
and hence 
\[(\forall n\geq m_1)(F_n(\bar{\alpha}^1\rest n)=F_n(\bar{\alpha}^2\rest
n)),\] 
finishing the proof of the claim. $\QED_{\ref{cl7}}$

\begin{claim}
\label{cl8}
If $\bar{\alpha}^1,\bar{\alpha}^2\in\pre\omega\lambda$ and $(\forall^\infty
n)(F_n(\bar{\alpha}^1\rest n)=F_n(\bar{\alpha}^2\rest n))$\\
then $(\forall^\infty n)(\alpha^1_n=\alpha^2_n)$
\end{claim}

\noindent{\em Proof of the claim:}\hspace{0.15in} Take $n_0<\omega$ such that
\[u(\bar{\alpha}^1)\cup u(\bar{\alpha}^2)\subseteq n_0\ \mbox{ and }\ (\forall
n\geq n_0)(F_n(\bar{\alpha}^1\rest n)= F_n(\bar{\alpha}^2\rest n)).\]
Then for all $n\geq n_0$ we have (by clause (a)  of the definition of
$F^0_{n_0,n_1,n_2}$):
\[\bl(\bar{\alpha}^1\rest n)=\bl(\bar{\alpha}^2\rest n)\quad\&\quad
\bk(\bar{\alpha}^1\rest n)=\bk(\bar{\alpha}^2\rest n).\]
Further, let $n_1>n_0$ be such that for all $n\geq n_1$,
$\bk(\bar{\alpha}^1\rest n)>n_0$. 

We are going to show that $\alpha^1_n=\alpha^2_n$ for all $n>n_1$. Assume not.
Then we have $n>n_1$ with $\alpha^1_n\neq \alpha^2_n$ and thus
$\eta_{\alpha^1_n}\neq\eta_{\alpha^2_n}$. Take $n'>n$ such that 
$\eta_{\alpha^1_n}\rest n'\neq \eta_{\alpha^2_n}\rest n'$. Applying \ref{cl4}
(2) and (4) choose $n''>n'$ such that 
\[\bm(\bar{\alpha}^1\rest n'')>n'\ \mbox{ and }\ k_0(\bar{\alpha}^1\rest
n'')=\max u(\bar{\alpha}^1).\] 
Now define inductively: $m_0=n''$, $m_{k+1}=\bm(\bar{\alpha}^1\rest m_k)$.\\
Thus
\[n''=m_0>\bl(\bar{\alpha}^1\rest m_0)\geq m_1>\bl(\bar{\alpha}^1\rest
m_1)\geq m_2>\ldots\] 
and 
\[m_k>\max u(\bar{\alpha}^1)\ \ \ \Rightarrow\ \ \ k_0(\bar{\alpha}^1\rest
m_k)=\max u(\bar{\alpha}^1)\] 
(see the definition of $\bm$). Let $k^*$ be the first such that $n\geq
m_{k^*}$ (so $k^*\geq 2$). Note that by the choice of $n_1$ above we
necessarily have
\[m_{k^*}>\bl(\bar{\alpha}^1\rest m_{k^*})=\bk(\bar{\alpha}^1\rest
m_{k^*-1})>n_0.\] 
Hence for all $k<k^*$:
\[\begin{array}{l}
F_{m_k}(\bar{\alpha}^1\rest m_k)=F_{m_k}(\bar{\alpha}^2\rest
m_k)\quad\mbox{ and }\\
\bl(\bar{\alpha}^1\rest m_{k+1})=\bl(\bar{\alpha}^2\rest m_{k+1})=
\bk(\bar{\alpha}^1\rest m_k)=\bk(\bar{\alpha}^2\rest m_k).\\
\end{array}\] 
Now, by the definition of the functions $\bl,\bm,\bk$ and the choice of $m_0$
(remember $k_0(\bar{\alpha}^1\rest m_0)=\max u(\bar{\alpha}^1)$) we know that
for each $i\in [\bk(\bar{\alpha}^1\rest m_k),\bl(\bar{\alpha}^1\rest m_k))$,
$k<k^*$ for some $\tau^m_\ell\in T_{\bl(\bar{\alpha}^1\rest m_k),m_k}$ and
$i_0,\ldots,i_{m-1}\in [\bl(\bar{\alpha}^1\rest m_k), m_k)$ we have
$\alpha^1_i=\tau^m_\ell(\alpha_{i_0}^1,\ldots,\alpha^1_{i_{m-1}})$. Moreover
we may demand that $\tau^m_\ell$ is a composition of depth at most
$\bl(\bar{\alpha}^1\rest m_k)-i$ of simple case terms. Since 
\[F^0_{\bk(\bar{\alpha}^1\rest m_k),\bl(\bar{\alpha}^1\rest
m_k),m_k}(\bar{\alpha}^1\rest m_k)= F^0_{\bk(\bar{\alpha}^2 \rest
m_k),\bl(\bar{\alpha}^2\rest m_k),m_k} (\bar{\alpha}^2\rest m_k)\]
we conclude that (by clause (d) of the definition of the functions
$F^0_{n_0,n_1,n_2}$): 
\[\alpha^2_i=\tau^m_\ell(\alpha^2_{i_0},\ldots,\alpha^2_{i_{m-1}}).\]
Now look at our $n$.\\
If $\bl(\bar{\alpha}^1\rest m_{k^*-1})>n$ then $\bk(\bar{\alpha}^1\rest
m_{k^*-1})\leq n<\bl(\bar{\alpha}^1\rest m_{k^*-1})$ and thus we find
$i_0,\ldots,i_{m-1}\in [\bl(\bar{\alpha}^1\rest m_{k^*-1}), m_{k^*-1})$ and
$\tau^m_\ell\in T_{\bl(\bar{\alpha}^1\rest m_{k^*-1}),m_{k^*-1}}$ such that
\[\alpha^1_n =\tau^m_\ell(\alpha_{i_0}^1,\ldots,\alpha_{m-1}^1)\ \ \&\ \
\alpha^2_n =\tau^m_\ell(\alpha_{i_0}^2,\ldots,\alpha_{m-1}^2).\] 
If $\bl(\bar{\alpha}^1\rest m_{k^*-1})\leq n$ then $n\in
[\bk(\bar{\alpha}^1\rest m_{k^*-2}), \bl(\bar{\alpha}^1\rest m_{k^*-2}))$
(remember that $\bl(\bar{\alpha}^1\rest m_{k^*-1})=\bk(\bar{\alpha}^1 \rest
m_{k^*-2})$ and $n<m_{k^*-1}\leq \bl(\bar{\alpha}^1\rest m_{k^*-2})$).
Consequently for some $i_0,\ldots,i_{m-1}\in [\bl(\bar{\alpha}^1\rest
m_{k^*-2}), m_{k^*-2})$ and  $\tau^m_\ell\in T_{\bl(\bar{\alpha}^1\rest
m_{k^*-2}),m_{k^*-2}}$ we have  
\[\alpha^1_n =\tau^m_\ell(\alpha_{i_0}^1,\ldots,\alpha_{m-1}^1)\ \ \&\ \
\alpha^2_n =\tau^m_\ell(\alpha_{i_0}^2,\ldots,\alpha_{m-1}^2).\] 
In both cases we may additionally demand that the respective term
$\tau^m_\ell$ is a composition of depth $\bl(\bar{\alpha}^1\rest m_{k^*-1})-n$
(or $\bl(\bar{\alpha}^1\rest m_{k^*-2})-n$, respectively) of terms of the
simple case. Now we proceed inductively (taking care of the depth of involved
terms) and we find a term $\tau\in T_{\bl(\bar{\alpha}^1\rest m_0),m_0}$
(which is a composition of depth at most $\bl(\bar{\alpha}^1\rest m_0)-n$ of
terms of the simple case) and $i_0,\ldots,i_{m-1}\in [\bl(\bar{\alpha}^1\rest
m_0),m_0)$ such that 
\[\alpha^1_n =\tau(\alpha_{i_0}^1,\ldots,\alpha_{m-1}^1)\ \ \&\ \
\alpha^2_n =\tau(\alpha_{i_0}^2,\ldots,\alpha_{m-1}^2).\] 
But now applying the clause (c) of the definition of the functions
$F^0_{n_0,n_1,n_2}$ we conclude that $\eta_{\alpha^1_n}\rest m_0
=\eta_{\alpha^2_n}\rest m_0$. Contradiction to the choice of $n'$ and the fact
that $m_0>n'$. $\QED_{\ref{cl8}}$
\medskip

The last two claims finish the proof of the theorem. $\QED_{\ref{second}}$

\begin{remark}
\label{kappafun}
{\em
If the models $M$ have $\kappa<\lambda$ functions (so $\langle\tau^n_i(x_0,
\ldots,x_{n-1}): i<\kappa\rangle$ lists the $n$--place terms) we can prove
$\KL(\lambda,\kappa)$ and the proof is similar.
}
\end{remark}
\medskip

\[*\qquad *\qquad *\]

\begin{finrem}
\label{finrem}
{\em
1)\quad Now we phrase exactly what is needed to carry the proof of theorem
\ref{main} for $\lambda>\kappa$. It is:
\smallskip 

\noindent $(\boxtimes)$\qquad for every model $M$ with universe $\lambda$ and
Skolem functions and with countable vcabulary, we can find pairwise distinct
$\alpha_{n,\ell}<\lambda$ (for $n<\omega, \ell<\omega$) such that
\begin{description}
\item[$(\otimes)$] if $m_0<m_1<\omega$ and $\ell_i^\prime<\ell_i^{\prime
\prime}$ for $i<m_0$ and $\ell_i<\omega$ for $i\in [m_0,m_1)$

then the models
\[\begin{array}{l}
({\rm Sk}(\{\alpha_{i,\ell^\prime_i},\alpha_{i,\ell^{\prime\prime}_i}: i<m_0
\}\cup \{\alpha_{m_0,k_0},\alpha_{m_0,k_1}\}\cup\{\alpha_{i,\ell_i}:i\in (m_0,
m_1)\}),\\
\qquad\alpha_{0,\ell^\prime_0},\alpha_{0,\ell^{\prime\prime}_0},\alpha_{1,
\ell^\prime_1},\alpha_{1,\ell^{\prime\prime}_1},\ldots,\alpha_{m_0-1,
\ell^\prime_{m_0-1}},\alpha_{m_0-1,\ell^{\prime\prime}_{m_0-1}},\alpha_{m_0,
k_0},\\
\qquad\alpha_{m_0,k_1},\alpha_{m_0+1,\ell_{m_0+1}},\ldots,\alpha_{m_1-1,\ell_{
m_1-1}})
  \end{array}\]
and 
\[\begin{array}{l}
({\rm Sk}(\{\alpha_{i,\ell^\prime_i},\alpha_{i,\ell^{\prime\prime}_i}: i<m_0
\}\cup \{\alpha_{m_0,k_0},\alpha_{m_0,k_2}\}\cup\{\alpha_{i,\ell_i}:i\in (m_0,
m_1)\}),\\
\qquad\alpha_{0,\ell^\prime_0},\alpha_{0,\ell^{\prime\prime}_0},\alpha_{1,
\ell^\prime_1},\alpha_{1,\ell^{\prime\prime}_1},\ldots,\alpha_{m_0-1,
\ell^\prime_{m_0-1}},\alpha_{m_0-1,\ell^{\prime\prime}_{m_0-1}},\alpha_{m_0,
k_0},\\
\qquad\alpha_{m_0,k_2},\alpha_{m_0+1,\ell_{m_0+1}},\ldots,\alpha_{m_1-1,\ell_{
m_1-1}})
  \end{array}\]
are isomorphic and the isomorphism is the identity on their intersection and
they have the same intersection with $\kappa$. 
\end{description}
For more details and more related results we refer the reader to
\cite{Sh:F254}. 

\noindent 2)\quad Together with \ref{koepke}, \ref{kappafun} this gives a good
bound to the consistency strength of $\neg\KL(\lambda,\kappa)$.

\noindent 3)\quad What if we ask $F_n:\pre n\lambda\longrightarrow
\pre{\omega>}\kappa$ such that $F_n(\eta)\trianglelefteq F_{n+1}(\eta)$ and
$\eta\in\pre\omega\lambda\ \Rightarrow\ F(\eta)=\bigcup F_n(\eta\rest n)\in
\pre\omega\kappa$? No real change.
}
\end{finrem}

\bigskip

\bibliographystyle{literal-unsrt}
\bibliography{listb,lista,listf,listx}

\shlhetal
\end{document}